\input amstex
\documentstyle{amsppt}
\magnification=1200
\hcorrection{0.2in}
\vcorrection{-0.4in}


\define\op#1{\operatorname{#1}}
\NoRunningHeads
\NoBlackBoxes
\topmatter
\title A New Proof of the Gromov's Theorem on Almost Flat Manifolds
\endtitle
\date (Dedicated to Jeff Cheeger's 75th Birthday)
\enddate
\author Xiaochun Rong \footnote{A part of this work was done during the author's sabbatical leave in the fall of 2018
at Capital Normal University, which was partially supported by NSFC Grant 11821101, Beijing Natural
Science Foundation Z19003, and a research fund from Capital Normal University. \hfill{$\,$}}
\endauthor
\address Mathematics Department, Rutgers University New Brunswick, NJ 08903 USA
\endaddress
\email rong\@math.rutgers.edu
\endemail
\address Mathematics Department, Capital Normal University, Beijing,
P.R.C.4.
\endaddress
\abstract We will give a new proof for the Gromov's theorem on almost flat manifolds, which is an inductive proof on dimension.
\endabstract
\endtopmatter
\document

The purpose of this paper is to give a new proof (by induction on dimension) for the well-known Gromov's theorem on almost flat manifolds.

\proclaim{Theorem 1} {\rm (\cite{Gr}, \cite{Ru})} Given $n\ge 2$, there are constants, $\epsilon(n), w(n)>0$,
such that if the sectional curvature and diameter of a compact $n$-manifold $M$ satisfies
$$|\op{sec}_M|\cdot \op{diam}(M)^2<\epsilon(n),$$
then $M$ is diffeomorphic to an infra-nilmanifold, $N/\Gamma$, where $N$ is a simply connected nilpotent Lie
group, $\Gamma$ is a discrete subgroup of $N\rtimes \op{Aut}(N)$ (the group of automorphisms)
such that the index $[\Gamma: \Gamma\cap N]\le w(n)$.
\endproclaim

A strong converse of Theorem 1 holds that any compact infra-nilmanifold admits one-parameter family of left invariant metrics (constructed via an inhomogeneous rescaling on a left invariant metric, determined by the infra-nil structure), $g_\epsilon$, such that $|\op{sec}_{g_\epsilon}|\cdot \op{diam}(g_\epsilon)^2\to 0$ as $\epsilon\to 0$ (\cite{Gr}, \cite{FH}).

Theorem 1 has been a cornerstone in the collapsing theory of Cheeger-Fukaya-Gromov on collapsed manifolds with bounded
sectional curvature (\cite{CFG}, \cite{CG1,2}, \cite{Fu1-3}), which has important applications (\cite{Fu4}, \cite{Ro1} and
references within).

In \cite{Gr}, Gromov proved that a bounded normal covering of $M$ is diffeomorphic to a nilmanifold, and that $M$
is diffeomorphic to an infra-nilmanifold was due to \cite{Ru} via constructing
a flat connection with a parallel torsion on $M$.

Roughly speaking, Gromov's proof in \cite{Gr} imitates the proof of Bieberbach theorem on a compact flat $n$-manifold $M$,
demonstrating that any deck transformation on $\Bbb R^n$ satisfies that its rotation component is either trivial or not small i.e., the
minimal rotation angle is bounded below by a positive constant depending only on $n$.
The core piece in Gromov's proof is a complicated distortion estimate on iterated commutators of deck transformations of short geodesic loops on the Riemannian universal covering $\tilde M$ that can have small non-trivial holonomy. The key discovery of Gromov is that when normalizing
the length of a short geodesic loop to one,  its holonomy is actually much smaller.

Gromov's proof involves several different fields and some unconventional arguments, which partially motivated the authors in \cite{BK}.

Our new proof of Theorem 1 is by induction on $n$: we show that to conclude
that an almost flat manifold has a bounded normal covering space that is diffeomorphic to a nilmanifold, it
is enough to explore the central part of a nilmanifold (see Theorem 3), provided with a recent criterion of
a nilmanifold (see Theorem 5). In particular, the induction enables us to bypasses the distortion estimate for iterated commutators of short geodesic loops (\cite{Gr}); indeed our proof reduces the above distortion estimate to a single special commutator; one of the short geodesic loops has an almost minimal length i.e., the length is proportional to the injectivity radius of $M$ (see (4.1)).

In our proof of Theorem 1, for simplicity we will always assume that an almost flat metric satisfies the following extra regularity (Theorem 2, also (4.2) below); an proof
without the extra regularity can be carried with a little more work (cf. \cite{Ro3}).

\proclaim{Theorem 2} {\rm (Smoothing)} Let $(M,g)$ be a complete $n$-manifold with $|\op{sec}_M|\le 1$.
For any $\epsilon>0$,  there is an $\epsilon$ $C^1$-close metric $g_\epsilon$ of $(\epsilon,l)$-regularity i.e.,
$$|g-g_\epsilon|_{C^1(M)}<\epsilon, \quad |\op{sec}_{g_\epsilon}|\le 1, \quad
|\nabla_\epsilon^k\op{Rm}(g_\epsilon)|\le c(n,k,\epsilon),\quad 0\le k\le l,$$
where $\op{Rm}$ denotes the curvature tensor, and $c(n,\epsilon,l)$ is a constant depending only on $n ,\epsilon$ and $l$.
\endproclaim

A proof of Theorem 2 via a Ricci flows can be found in \cite{Shi}, or via a local embedding method was given in \cite{Ab}.

The following theorem is the key in the proof of Theorem 1 via induction on $n$.

\proclaim{Theorem 3} {\rm (Central geometry/topology of almost flat manifolds)} There exists a constant $\epsilon(n)>0$ such that if $M$ is an almost flat $n$-manifold i.e.,
$|\op{sec}_M|\cdot (\op{diam}(M))^2<\epsilon\le \epsilon(n)$, then the following properties hold:

\noindent {\rm (3.1)} A finite normal covering space of $M$ with order $\le a(n)$ admits a principal $T^k$-bundle, $T^k\to \hat M\to \hat M/T^k$, and a $C^1$-close $T^k$-invariant metric such that the quotient metric on $\hat M/T^k$ is $\epsilon'$-almost flat, $0<\epsilon'\le \epsilon(n)$.

\noindent  {\rm (3.2)} Normalizing the metric to that $|\op{sec}_M|\le 1$ (thus $\op{diam}(M)^2<\epsilon(n)$), the injectivity
radius of the Riemannian universal cover, $\op{injrad}(\tilde M)\ge \delta(n)>0$.
If in addition, $\pi_1(M)$ is nilpotent, then a short basis has almost constant
displacement everywhere in a fundamental domain.
\endproclaim

\remark{\rm Remark 4} (4.1) The proof of (3.1) requires a distortion estimate only for a special commutator of two short geodesic loops; one of which has length proportional to $\op{injrad}(M)$ by a constant (see (8.1)).

\noindent (4.2) The $T^k$-action on $\hat M$ almost preserves the almost flat metric; averaging which under the $T^k$-action yields a desired invariant metric in (3.1) (here a $(\epsilon,l)$-regularity in Theorem 2 is used, say $l=2$). We point it out that without 
Theorem 2, an inductive proof can still be carried out, with additional work (\cite{Ro3}),
\endremark

Theorem 3 enables us to prove Theorem 1 by induction on $n$, with a help from the following topological
criterion for a compact nilmanifold.

\proclaim{Theorem 5} {\rm (Nilmanifolds: Iterated principal circle bundles, \cite{Na}, \cite{Be})}
A compact $n$-manifold $M$ is diffeomorphic to a nilmanifold if and only if $M$ admits an iterated principal circle bundles,
$$S^1\to M\to M_1, \quad S^1\to M_1\to M_2,\quad \cdots, \quad S^1\to M_n\to \op{pt}.$$
\endproclaim

By induction on $n$, it is easy to show that a nilmanifold admits an iterated principal circle bundles, and conversely,
that the fundamental group is nilpotent and the universal cover is diffeomorphic to $\Bbb R^n$ (hence $M$ is
homeomorphic to a nilmanifold, \cite{FH}).  A diffeomorphism in Theorem 5 was mentioned in early literature
without a proof (cf. \cite{Be}).  The author observed that the main result in \cite{Na} implies Theorem 5, and an elementary proof can be found in \cite{Be}.

Let's first give a proof of Theorem 1 by assuming Theorem 3.

\demo{Proof of Theorem 1}

The proof is divided into two steps: Step 1. Prove that a bounded normal covering space of $M$ is diffeomorphic to a nilmanifold.
Step 2. Prove that $M$ is diffeomorphic to an infra-nilmanifold.

Step 1. We proceed by induction on $n$, starting with the trivial case: $n=2$.

For $n>2$, by Theorem 3 we conclude that $M$ has a bounded normal covering space $\hat M$ which
admits a principal $T^k$-bundle ($k\ge 1$), $T^k\to \hat M@>f>>\hat M/T^k=B$, and a $C^1$-close $T^k$-invariant metric, $g_{T^k}$, such that the quotient metric on $B$ is almost flat. Without loss of generality, we may assume $k<n$.
Applying induction on $B$, we may assume that $B$ has a bounded normal covering space, $\hat \pi: \hat B\to B$,
and $\hat B$ is diffeomorphic to a nilmanifold. Let  $T^k\to \hat \pi^*(\hat M)@> >>\hat B$ denote the
$\hat \pi$-pullback principal $T^k$-bundle.  Then  $M':=\hat \pi^*(\hat M)$ is a bounded normal covering space of $M$.

By Theorem 5, $\hat B$ admits iterated principal circle bundles, thus $M'$ admits an iterated principal circle bundles.
Again by Theorem 5, $M'$ is diffeomorphic to a nilmanifold, $N/\Gamma'$, where $N$ is a simply connected nilpotent
manifold which contains $\Gamma'=\pi_1(M')$ as a co-compact discrete subgroup of $N$.

\vskip2mm

Step 2.  Based on the above Step 1, we assume a $\Gamma'$-conjugate diffeomorphism, $\phi: (\tilde M,\Gamma')\to (N,\Gamma')$. Because
$\Gamma$ is torsion free, by Malc\'ev rigidity the homomorphism induced by conjugation, $\Gamma\to \op{Aut}(\Gamma')$, embeds
$\Gamma$ into $N\rtimes \op{Aut}(\Gamma')$, hence $N/\Gamma$ is an infra-nilmanifold.  We shall extend $\phi$, using
the center of mass method, to a $\Gamma$-conjugate diffeomorphism, $(\tilde M,\Gamma)\to (N,\Gamma)$, thus $M$ is diffeomorphic
to $N/\Gamma$.

By the above we may assume the following: identifying $M'$ with $N/\Gamma'$, $M'$ admits two finite group $\Gamma/\Gamma'$-actions: one isometric free $\Gamma/\Gamma'$-action
that almost preserves the nilpotent structure i.e., the iterated principal $S^1$-bundles,
while the other preserves the nilpotent structure. 
To apply the center of mass, it is enough to construct a left invariant metric on $N/\Gamma'$ such that $\Gamma/\Gamma'$ acts on $N/\Gamma'$ isometrically and the two metrics are $C^1$-close.
A standard method is averaging the almost flat metric (of higher regularity) on $\tilde M$ by the nilpotent Lie group action to obtain a left invariant metric.

Based on (3.2) one may give a direct construction of a left invariant metric. Observe that the Gromov's short basis at $p'\in M'$, $\{\gamma_j\}$, determines a basis on the Lie algebra $\hbar=T_{\tilde p}\tilde M$, $\{v_k\}$, $v_k=\gamma_k'(0)$ with $|\gamma_k'(0)|=\op{length}(\gamma_k)$.
We define a left-invariant metric on $N$ by
$$\tilde g_N(v_k,v_l)=\tilde g(v_k,v_l).$$
By abusing notation, we use $g_N$ to denote the averaging $g_N$ by $\Gamma/\Gamma'$.
By (3.2), it is clear that restricting to any $R$-unit ball
on $N$ with respect to $\tilde g_N$ or $\tilde g$, the two metrics are GH-close to a flat
metric, thus the two metrics are GH-close on $B_R(\tilde p)$.

Using the two $C^1$-close metrics on $M'$, $g'$ and $g_N$, and identifying $\tilde M$ with $N$ by $\phi$, we may assume
two free isometric $\Gamma$-actions on  $\tilde M$, denoted by $\mu_1$ and $\mu_2$ such that $\mu_1\equiv \mu_2$ on $\Gamma'$, with respect to $\tilde g$ and $\tilde g_N$ respectively. To apply the method of the center of mass, we normalize $\tilde g$ (and $\tilde g_N$) to that $|\op{sec_{\tilde g}}|\le 1$, thus $\op{diam}(M')^2<a(n)^2\epsilon(n)<\frac {\delta(n)}2$, the convexity radius of $\tilde g$ (see (3.2)).

Let $\Gamma=\bigcup_{k=1}^s \alpha_k\Gamma'$, $s\le a(n)$. For fixed $\alpha=\alpha_k$, $h\in N$, and any $\gamma\in \Gamma'$,
$$\split d(\mu_2(\alpha\gamma)(\mu_1((\alpha\gamma)^{-1},h),h)&=d(\mu_2(\alpha)\mu_2(\gamma)   (\mu_1(\gamma^{-1})\mu_1(\alpha^{-1})h)),h)\\&=d(\mu_2(\alpha)\mu_1(\alpha^{-1})h,h).\endsplit$$
Then $A(h)=\{\mu_2(\alpha_k\gamma)^{-1}\mu_1(\alpha_k\gamma)(h),\,\gamma\in \Gamma', 1\le k\le s\}$ is a finite set of size $|A(h)|\le a(n)$.
We now specify a representative $\alpha_i$ as a projection of the identity $e$ to $\alpha_i\Gamma'$ (with respect to $\tilde g$). Because $\op{diam}(M')<c(n)\epsilon(n)$ is small,
$$\delta(n)>10a(n)\max\{d(\mu_2(\alpha_k)\mu_1(\alpha^{-1}_k)(h,h), \, 1\le k\le s\},$$
thus we can define a map, $f: N\to N$, $f(h)=c_h$, where $c_h$ denotes the center of mass of $A(h)\subset B_s(h)$ with respect to $\tilde g$.
Clearly, $f$ is $\Gamma$-invariant,  and it is straightforward to check that for any $0\le t\le 1$, the map, $F(t,h)=c_h(t): [0,1]\times \tilde N\to N$,  defines an isotopic between $F(0,\cdot)=\op{id}_N$ and $F(1,\cdot)=f$, where $c_h(t)$ is the unique minimal geodesic from $h$ to $f(h)$ (cf. \cite{GK}).
\qed\enddemo

In the rest of the paper, our main effort is to prove Theorem 3.
Because $|\op{sec}_M|\cdot \op{diam}(M)^2$ is a scaling invariant, by a standard Gromov's compactness argument Theorem 3 is equivalent to the following:

\proclaim{Theorem 3'} {\rm (Central geometry/topology of almost flat manifolds)}  Let a sequence of compact $n$-manifolds, $M_i@>\op{GH}>>\op{pt}$,
such that $|\op{sec}_{M_i}|\le 1$. Passing to a subsequence if necessary, the following properties hold:

\noindent {\rm (3.1)'} $M_i$ has a bounded normal covering space $\hat M_i$
which admits a principal $T^k$-bundle  ($k\ge 1$) and a $C^1$-close invariant metric so that the base manifolds,
equipped with the quotient metrics, satisfy that
$\hat M_i/T^k@>\op{GH}>>\op{pt}$ and $|\op{sec}_{\hat M_i/T^k}|\le c(n)$, a constant depends on $n$.

\noindent {\rm (3.2)'} The injectivity radius of the Riemannian universal covering space of $M_i$, $\op{injrad}(\tilde M_i)\ge \delta(n)>0$. If in addition, $\pi_1(M_i)$ is nilpotent, then a short basis at any $p_i\in M_i$, the ratio of displacement, $\frac{d(\gamma_i(\tilde p_i),\tilde p_i)}{d(\gamma_i(\tilde x_i),\tilde x_i)}\to 1$, for all $\tilde x_i$ in
a fundamental domain at $\tilde p_i\in \tilde M_i$.
\endproclaim

In our proof of (3.1)', we blow-up twice on the sequence, $M_i$, by $\ell_i^{-1}$ and $\rho_i^{-1}$, where $\ell_i=\op{diam}(M_i)$ (the minimal collapsing rate) and $\rho_i=\op{in jrad}(M_i)$ (the maximal collapsing rate) respectively, and we investigate structures on the blow-up
sequences by standard tools in metric Riemannian geometry (e.g., the equivariant GH-convergence, the Cheeger-Gromov convergence theorem, gluing of $C^1$-close local fiber bundles via the method of the center of mass, etc). Structural results from the blow-up by
$\ell_i^{-1}$ will be used in exploring on a finite normal cover of $\rho_i^{-1}M_i$, a local product structure of a $T^k$-bundle (see Example 10).

Because the proof is fairly involved, for the convenience of readers we first present an
outline: let's start with the following $\ell_i^{-1}$-blow up commutative equivariant GH-convergence,
$$\CD (B_{\frac \pi{\ell_i}}(0_i),0_i, \Gamma_i)@>\op{eqGH}>>(\Bbb R^n,0,G)\\
@V\exp_{p_i}VV @V\op{proj}VV\\
(\ell_i^{-1}M_i,p_i)@>\op{GH}>>X=\Bbb R^n/G,\endCD$$
where via path lifting $\Gamma_i=\pi_1(M_i,p_i)$ acts pseudolly on $B_{\frac \pi{\ell_i}}(0_i)\subset T_{p_i}M_i$. Because
$\op{diam}(X)=1$ (compact), using the generalized Bieberbach theorem (see Theorem 6 below),
by the standard equivariant convergence it follows that a bounded normal covering space of $M_i$, $\ell_i^{-1}\hat M_i@>\op{GH}>>T^m$, a flat torus. Consequently, one constructs a fiber bundle, $N_i\to \ell_i^{-1}\hat M_i@>f_i>>T^m$, such that $f_i$ is an
$\epsilon_i$-Gromov-Hausdorff approximation (briefly, GHA) and $\epsilon_i$-Riemannian submersion, and the second fundamental form $|\op{II}(N_i)|\le c(n)$ (with respect to the original metric), see Theorem 7 (a new proof in Appendix) and Lemma 8 below (cf. \cite{Fu1}, \cite{CFG}).

Because an $f_i$-fiber is almost flat, we establish, by induction on $n$, the following properties ($i$ large):  i) $\tilde M_i$ is diffeomorphic to $\Bbb R^n$, and (3.2) holds. ii) any $x_i\in \Hat M_i$, $c(n)^{-1}\le \op{injrad}(\hat M_i)/\op{injrad}(x_i)\le c(n)$.  iii) Based on a holonomy distortion estimate for short geodesic loops whose lengths are proportional to $\rho_i$ (see Lemma 9), we
can choose $\hat M_i$ so that the pointed GH-limit,  $(\rho_i^{-1}\hat M_i,x_i)@>\op{GH}>>(T^k\times \Bbb R^{n-k}, e\times 0)$,
a product of flat $n$-manifolds, such that $k$ is independent of $x_i$, $\op{diam}(T^k)\le d$ and $\op{injrad}(T^k)\ge d^{-1}$
($d$ is a constant independent of $x_i$).   Applying the Cheeger-Gromov's convergence  theorem (\cite{Ch}, \cite{GLP}),
we construct a locally finite open cover
for $\hat M_i$, $\{U_{i,\alpha}\}$, such that $\rho_i^{-1}U_{i,\alpha}$ is $C^3$ $\epsilon_i$-close to a product of flat manifolds,
$T^k\times \b B_{10}^{n-k}(0_\alpha)$.

Observe that a choice of a base point $e$ on the $T^k$-factor of $T^k\times \b B^{n-k}_{10}(0_\alpha)$ and
a canonical basis for $\pi_1(T^k,e)$ uniquely determines a group structure on $T^k$-fibers (via parallel translations),
hence $T^k\times \b B^{n-k}_{10}(0_\alpha)$ is a principal $T^k$-bundle. Moreover, a different choice of $e$ and a canonical
basis for $\pi_1(T^k,e)$ yields an automorphism of $T^k$, thus an isomorphism of the resulting principal $T^k$-bundles.
If $U_{i,\alpha}\cap U_{i,\beta}\ne \emptyset$,  it is clear that the two trivial $T^k$-bundles are $C^1$-close, thus we
are able to glue $\{U_{i,\alpha}\}$ to form a $T^k$-bundle, $T^k\to \hat M_i\to B_i$, with an affine structural group
(note that our gluing operation can be viewed as a simplest case of the construction of an T-structure on collapsed
manifolds with bounded sectional curvature, \cite{CG2}).  Because $\pi_1(T^k)$ injects into the center of $\pi_1(\hat M_i)$,
this is a principal $T^k$-bundle.

We begin to fill in details in the above outline of proof of Theorem 3'.

\proclaim{Theorem 6} {\rm (Generalized Bieberbach Theorem, \cite{FY})}
Let $G$ be a closed subgroup of $\op{Isom}(\Bbb R^n)$ such that $\Bbb R^n/G\ne \{\op{pt}\}$ is compact. Then

\noindent {\rm (6.1)} The identity component $G_0$ of $G$ satisfies that $\Bbb R^n/G_0$ with the quotient metric is
isometric to $\Bbb R^m$, $m=n-\dim(G_0)$.

\noindent {\rm (6.2)} $G$ contains a normal subgroup $\hat G$ of finite index such that $\Bbb R^n/\hat G$ is isometric to
a flat torus $T^m$.
\endproclaim

We present a short proof here, based on basic properties of a complete non-compact length space with non-negative curvature in the sense of Toponogov triangle comparison (\cite{BGP}), which are analogous to the classical Soul theorem and Splitting theorem of Cheeger-Gromoll in Riemannian geometry (\cite{ChG1,2}).

\demo{Proof of Theorem 6}

(6.1) Arguing by contradiction, assume that $\Bbb R^n/G_0$ is not isometric to $\Bbb R^m$ ($1\le m<n$). Then $\Bbb R^n/G_0$, as a complete non-compact Alexandrov space of non-negative curvature, splits, $\Bbb R^n/G_0=\Bbb R^k\times X$ ($0\le k<m$), and $X$ contains no line (\cite{Mi}). Via horizontal lifting of lines in $\Bbb R^k$, $\Bbb R^n$ splits as, $\Bbb R^n=\Bbb R^k\times \Bbb R^{n-k}$, such that each $\Bbb R^{n-k}$-factor is $G_0$-invariant. We claim that if $X=\Bbb R^{n-k}/G_0$ is not a point, then $X$ is not compact. Without loss of generality,
we may assume that the $G_0$-orbits do not form a fiber bundle i.e., there is a non-principal
$G_0$-orbit. If $p$ is a point with a local maximal isotropy group, then
$G_0(p)$ is a subplane $\Bbb R^\ell$, $1\le \ell< n-k$. Without loss of the generality we may assume $p=0$. Let $H=\{A \in SO(n), \, (A,b)\in G_0\}$, and $\bar H$ denote the closure of $H$ in $SO(n)$. Let $G_0(0)=\Bbb R^k$ ($k\ge 1$,
which includes the case that $\Bbb R^k\rtimes SO(k)<G_0$), and let $\Bbb R^{n-k}$ denote the orthogonal complement of $\Bbb R^k$. Then for any $z\in S^{n-1}_1$, 
$G_0(z)\subset \Bbb R^k\times S^{n-k-1}$. Consequently for
$G_0(z_i)\subset \Bbb R^k\times S^{n-k-1}_{R_i}$, $d(G_0(z_1),G_0(z_2))\ge |R_1-R_2|$, thus the $G_0$-orbit space is not compact. 

Because $X$ is not compact and contains no line, $X$ has a soul $S$. Let $\gamma$ be a ray in $X$ from $p\in S$. Because for any $\alpha\in G/G_0$, $\alpha(p)\in P$ and $\alpha(\gamma)$ is a ray in $X$, $\Bbb R^n/G=\Bbb R^m/(G/G_0)$ contains a ray, a
contradiction to that $\Bbb R^m/(G/G_0)=\Bbb R^n/G$ is compact.

(6.2) Because $\Bbb R^n/G$ is compact, $G/G_0$ is finitely generated. By (6.1) and 6.13 in \cite{Ra}, $G/G_0$ contains a torsion-free
subgroup $\Gamma$ of finite index.  Without loss of the generality, $\Gamma$ can be taken to be normal in $G/G_0$,
thus $\Bbb R^m/\Gamma$ is a compact flat manifold. Now we apply the Bieberbach theorem to obtain a normal subgroup $\Bbb Z^m$
of $\Gamma$, with the desired normal subgroup of $G$ as the inverse image of the normalizer of $\Bbb Z^m$ in $G/G_0$.
\qed\enddemo

\proclaim{Theorem 7} {\rm (Fiber bundles, \cite{Fu1}, \cite{CFG})} Given $n, d, \rho>0$, there exist constants, $\epsilon(n,d,\rho), c(n)>0$, such that if a compact $n$-manifold $M$ and
a compact $m$-manifold $N$ satisfy the following conditions:
$$|\op{sec}_M|\le 1,\quad |\op{sec}_N|\le 1,\quad \op{diam}(M)\le d,\quad \op{injrad}(N)\ge 2\rho,$$
and $d_{\op{GH}}(M,N)<\epsilon\le \epsilon(n,d,\rho)$, then there is smooth fiber bundle map,
$F\to M@>f>>N$, satisfying the following properties:

\noindent {\rm (7.1)} $f$ is an $\epsilon$-GHA.

\noindent {\rm (7.2)} $f$ is an $\epsilon'$-Riemannian submersion i.e., for any vector $\xi$ orthogonal to $F$,
$e^{-\epsilon'}\le \frac{|df(\xi)|}{|\xi |}\le e^{\epsilon'}$, where $\epsilon'=\Psi(\epsilon|n,d,\rho)\to 0$ as
$\epsilon\to 0$.

\noindent {\rm (7.3)} The second fundamental form of any fiber, $|\op{II}(F)|\le c(n)$.
\endproclaim

Note that by (7.3), the intrinsic metric on $F$ is almost flat, and the structural group can
be reduced to an affine group (\cite{Fu1}, \cite{CFG}).

There have been various extension of Theorem 7 under weak curvature conditions; see \cite{Ya}
for $\op{sec}_M\ge -1$, and see \cite{Hu} and \cite{Ro3} for $\op{Ric}_M\ge -(n-1)$ and local
covering space non-collapsed. In Appendix, we will present a self-contained elementary proof.

In the next two lemmas, we will establish properties that will be used in the proof of Theorem 3'.

\proclaim{Lemma 8} Let the assumption be as in Theorem 3'.  Passing to a subsequence if necessary, the following properties hold:

\noindent {\rm (8.1)} There is a bounded normal covering space of $M_i$, $\ell_i^{-1}\hat M_i@>\op{GH}>> T^m$, a flat torus.

\noindent {\rm (8.2)}  For $i$ large, there is a smooth fiber bundle, $N_i\to \hat M_i@>f_i>> \ell_iT^m$, satisfying {\rm (7.1)}-{\rm (7.3)}
(with respect to the original metric).

\noindent {\rm{(8.3)} The universal cover of $\tilde M_i$ is diffeomorphic to $\Bbb R^n$.

\noindent{\rm (8.4)}  The injectivity radius, $\op{injrad}(\tilde M_i)\ge \delta(n)>0$.

\noindent {\rm (8.5)} For any $x_i \in M_i$, $1\le \frac{\op{injrad}(x_i,M_i)}{\op{injrad}(M_i)}\le b(n)$.

\noindent {\rm (8.6)} For any $x_i\in M_i$, $(\rho_i^{-1}M_i,x_i)@>\op{GH}>>(X,x)$. Then $X$ has a soul $F$ of dimension $k\ge 1$,
$\op{diam}(F)\le d$ and $\op{injrad}(F)\ge d^{-1}$, where $k$ and $d$ are constants independent of $x_i$ (which may depend on
the sequence $\{M_i\}$).

\noindent {\rm (8.7)} $\hat \Gamma_i=\pi_1(\hat M_i)$ contains a subgroup $\Bbb Z^k$, generated by short geodesic loops of length uniformly proportional to $\rho_i$, with normalizer of index, $[\hat \Gamma_i: N_{\hat \Gamma_i}(\Bbb Z^k)]\le c(n)$.
\endproclaim

\demo{Proof} (8.1) By the equivariant Gromov's pre-compactness, passing to a subsequence we may assume the
following commutative diagram (see the outline of the proof of Theorem 3.2'):
$$\CD (B_{\frac \pi{\ell_i}}(0_i),0_i, \Gamma_i)@>\op{eqGH}>>(\Bbb R^n,0,G)\\
@V\exp_{p_i}VV @V\op{proj}VV\\
(\ell_i^{-1}M_i,p_i)@>\op{GH}>>X=\Bbb R^n/G.\endCD$$

By Theorem 3.5, we may assume a normal subgroup of $G$ of finite index, $\hat G$, such that $\Bbb R^n/\hat G=T^m$ and
$G/\hat G$ is finite. Then $\Gamma_i$ contains a normal subgroup, $\hat \Gamma_i$,
such that $\hat \Gamma_i\to \hat G$, $\hat \Gamma_i/\hat \Gamma_i\cong G/\hat G$,
thus $\hat M_i=B_{\frac \pi{\ell_i}}(0_i)/\hat \Gamma_i$ is a normal covering space of $M_i$ of order $|G/\hat G|$,
and $\ell_i^{-1}\hat M_i@>\op{GH}>>T^m$ (let $\{\alpha_j\}$ be a finite set of generators for $G/G_0$, and let $\alpha_{i,j}\in \Gamma_i$ such that $\alpha_{i,j}\to \alpha_j$. Then
$\hat \Gamma_i=\left<\Gamma_i(\epsilon),\alpha_{i,j}\right>$, cf. \cite{FY}).

\vskip1mm

(8.2) can be obtained from (8.1) and Theorem 7; note that $d_{\op{GH}}(\ell_i^{-1}\hat M_i,T^m)<\epsilon_i$ implies that
$d_{\op{GH}}(\hat M_i,\ell_iT^m)<\ell_i\epsilon_i$, where $\op{sec}_{\ell_iT^m}\equiv 0$ and $\op{injrad}(\ell_iT^m)=\ell_i\op{injrad}(T^m)$.
Because
$$\frac{d_{\op{GH}}(\hat M_i,\ell_iT^m)}{\op{injrad}(\ell_iT^m)}<\frac {\epsilon_i\ell_i}{\ell_i\op{injrad}(T^m)}=\frac {\epsilon_i}{\op{injrad}(T^m)}\to 0,$$
Theorem 7 applies to the pair, $(\hat M_i,\ell_iT^m)$, with fixed $i$ large, to obtain a
bundle map, $f_i: \hat M_i\to \ell_i T^m$, which is the same obtain by rescaling by $\rho_i^{-1}$.
To see that $|\op{II}(N_i)|\le c(n)$ with respect to the original metric, equivalently $|\op{II}(N_i)|\le c(n)\ell_i\to 0$ with respect to $\ell_i^{-1}\hat M_i$, we will assume the higher regularity of the original metric (Theorem 3) that guarantees that $\hat b_{i,\beta,j}$ is $C^{2,\alpha}$-close to $\b b_j$ (see Remark A6).


\vskip1mm

In the proof of (8.3)-(8.7), we will proceed induction on $n$, using the fiber bundle structure obtained in (8.2).

(8.3) We proceed by induction on $n$, starting with the trivial case that $n=2$. For $n>2$, by (8.2)  we obtain a smooth fiber bundle map,
$N_i\to \hat M_i@>f_i>>\ell_iT^m$ and $f_i$ satisfying (7.1) and (7.2).  Let  $\pi: \Bbb R^m\to \ell_iT^m$ denote the Riemannian universal cover, and let $N_i\to \pi^*(\hat M_i)@>\hat f_i>>\Bbb R^m$ denote the pullback bundle by $\pi$.  Then $\pi^*(\hat M_i)$ is diffeomorphic to $N_i\times \Bbb R^m$, and $\pi^*(\hat M_i)$ is a covering space of $\hat M_i$.  Hence $\tilde M_i$ is diffeomorphic to $\tilde N_i\times \Bbb R^m$, where $\tilde N_i$ is a universal covering space of $N_i$. Because the intrinsic metric on $N_i$ satisfies that
$N_i@>\op{GH}>>\op{pt}$ such that $|\op{sec}_{N_i}|\le c(n)$, applying the induction on $N_i$ we conclude that $\tilde N_i$ is diffeomorphic
to $\Bbb R^{n-m}$, thus $\tilde M_i$ is diffeomorphic to $\Bbb R^n$.

\vskip1mm

(8.4) As seen from the proof of (8.3), $\tilde M_i$ has a fiber bundle, $\tilde N_i\to \tilde M_i@>\tilde f_i>>\Bbb R^m$, such that $\tilde f_i$ also satisfies (7.2) and (7.3).

We proceed by induction on $n$, starting with $n=2$, thus $m=1$ or $2$.  For $m=1$, the above implies that $\tilde M_i$ is a $\Bbb R^1$-bundle over $\Bbb R^1$. If $\tilde M_i$ has a short geodesic loop $\tilde \gamma_i$ at $\tilde x_i$, then $\tilde \gamma_i$ can be shortly homotopic to a trivial loop at $\tilde x_i$, a contradiction.  If $m=2$, then $\ell_i^{-1}\hat M_i$ is diffeomorphic and $C^{1,\alpha}$-close to a flat torus $T^2$, thus $g_i$ is bi-Lipschitz to the flat metric on $T^2$. By loop shortening operation on
the projection of $\tilde \gamma_i$ in $\hat M_i$ with respect to the flat metric, one gets a non-trivial geodesic loop with
respect to the flat metric, a contradiction because any geodesic loop on a flat $T^2$ is not homotopically trivial.

For $n>2$, by applying induction on $N_i$, we may assume that the intrinsic metric on $\tilde N_i$ satisfies $\op{injrad}(\tilde N_i)\ge \delta_1(n)$.  If $\tilde M_i$ has a short geodesic loop $\tilde \gamma$ at $\tilde x_i$, by (7.2) and (7.3) one sees that via the radial contraction in a $1$-tube around $\tilde N_i\ni \tilde x_i$, followed by a loop shortening process in $\tilde N_i$, one shortly deforms $\tilde \gamma_i$ to an intrinsic short geodesic loop in $\tilde N_i$, thus whose lifting in $\tilde N_i\subset \tilde M_i$ is not a loop, a contradiction.

\vskip1mm

(8.5) Observe that (8.5) holds if it holds on a bounded normal covering space of $M_i$ (i.e., $\hat M_i$ in (8.1)).

As a preparation, we first show that for any $x_i\in \hat M_i$, $1\le \frac {\op{injrad}(x_i,N_i)}{\op{injrad}(x_i,\hat M_i)}\le 2$. By definition,
$\op{injrad}(x_i,N_i)\ge \op{injrad}(x_i,\hat M_i)$. If $\gamma_i$ is a geodesic loop at $x_i$ with $|\gamma_i|=2\op{injrad}(x_i,\hat M_i)$,
then $\gamma_i$ is contained in a $r$-tubular neighborhood of $N_i$, $r<\op{diam}(M_i)$. By the radial deformation to $N_i$, followed by a path shortening one deforms $\gamma_i$ shortly to a $N_i$-geodesic loop at $x_i$ with almost the same length, because $|\op{II}(N_i)|\le c(n)$ and
$r<<1$ (the length distortion is at most exponentially in terms of $c(n)r$).

We now proceed by induction, starting with a trivial case $n=2$. As seen in (8.3) and (8.4), we consider the fiber bundle, $N_i\to \hat M_i@>f_i>>\ell_iT^m$. Applying induction on $N_i$, we assume that $N_i$ has a bounded normal covering space $\hat N_i$ such that $1\le \frac{\op{injrad}(x_i,\hat N_i)}{\op{injrad}(\hat N_i)}\le b_1(n)$.
Consequently, $1\le \frac{\op{injrad}(x_i,N_i)}{\op{injrad}(N_i)}\le b_2(n)$.

Observe that for any $p_i, x_i\in \hat M_i$, let $\alpha_i$ be a minimal geodesic from $f_i(p_i)$ to $f_i(x_i)$. Let $\tilde \alpha_i$ denote
the horizontal lifting of $\alpha_i$ at $p_i$, $q_i=\alpha_i(1)$ and $x_i$ are in the same $N_i$-fiber.
Because $|\op{II}(N_i)|\le c(n)$ and $d(f_i(p_i),f_i(x_i))<\op{diam}(M_i)<<1$, similarly we have that
$2^{-1}\le \frac{\op{injrad}(p_i,\hat M_i)}{\op{injrad}(\tilde \alpha_i(1),\hat M_i)}\le 2$.

Assume that $\rho_i=\op{injrad}(p_i,\hat M_i)$. For any $x_i\in \hat M_i$, let $q_i=\tilde \alpha_i(1)$.
Then
$$\split 1\le \frac {\op{injrad}(x_i,\hat M_i)}{\op{injrad}(p_i,\hat M_i)}&=\frac {\op{injrad}(x_i,\hat M_i)}{\op{injrad}(q_i,\hat M_i)}\cdot \frac{\op{injrad}(q_i,N_i)}{\op{injrad}(p_i,\hat M_i)}\cdot \frac {\op{injrad}(q_i,N_i)}{\op{injrad}(p_i,N_i)}\cdot \frac{\op{injrad}(p_i,N_i)}{\op{injrad}(p_i,\hat M_i)}\\&\le 2^4\frac{\op{injrad}(x_i,N_i)}{\op{injrad}(q_i,N_i)}\le 2^4b_2(n).
\endsplit$$

\vskip1mm

(8.6) Similar to (8.5), (8.6) holds if it holds on a sequence of bounded normal covering space of $M_i$, $\hat M_i$. Assume that $\rho_i=\op{injrad}(p_i,\hat M_i)$.

We proceed by induction on $n$, starting with the obvious case that $n=2$. For $n>2$, consider $N_i\to \hat M_i@>f_i>>\ell_iT^m$ with (7.1) and (7.2). Without loss of generality, we assume that $\frac{\rho_i}{\ell_i}\to 0$ (otherwise $M_i$ is diffeomorphic to a flat manifold); so $(\rho_i^{-1}\hat M_i,x_i)@>\op{GH}>>(X,x)$ with $X$ a non-compact flat $n$-manifold. Because $|\op{II}(N_i)|_{\rho_i^{-1}\hat M_i}\le c(n)\rho_i\to 0$,  $X=Y\times \Bbb R^m$ (Splitting theorem of Cheeger-Gromoll), and $(\rho_i^{-1}N_i,x_i)@>\op{GH}>>(Y,x)$.  By (8.5),  it is clear that $Y$ has a non-trivial geodesic
loop at $x$, thus a soul of $Y$, $F$, satisfies that $d(x,F)<\infty$ and $\dim(F)\ge 1$.  Applying induction on $N_i$, we obtain
that $\op{diam}(F)\le d$ and $\op{injrad}(F)\ge d^{-1}$.

\vskip1mm

(8.7) In the proof of (8.6) we obtain that $N_i\to \hat M_i\to \ell_iT^m$, $(\rho_i^{-1}\hat M_i,x_i)@>\op{GH}>>(Y\times \Bbb R^m,x)$ and $(\rho_i^{-1}N_i,x_i)@>\op{GH}>>(Y,x)$, $Y$ has a
soul $F$, a compact flat submanifold. By Cheeger-Gromov $C^{1,\alpha}$-convergence, we
may assume that $F_i\subset N_i$ and $B_{d(F,x)}(F_i)$ is $C^{1,\alpha}$-close to $B_{d(F,x)}(F)$ (in particular, $F_i$ is $C^{1,\alpha}$-close to $F$). By Bieberbach theorem, $\pi_1(F_i)\cong \pi_1(F)$ contains a normal subgroup, $H_i\cong\Bbb Z^k$, with $[\pi_1(F_i):H_i]=a\le a(k)$. We claim that $\pi_1(F_i)$ is normal in $\hat \Gamma_i$, assuming which we
first prove (3.6.7).

Let $A_i$ denote the conjugate class of $H_i$ in $\hat \Gamma_i$, and let $S(A_i)$ be the group of permutations on $A_i$. Let $\phi: \hat \Gamma_i\to S(A_i)$ denote the homomorphism induced by conjugation. Then $N_{\hat\Gamma_i}(H_i)\supseteq \ker \phi$, thus $[\hat \Gamma_i: N_{\hat \Gamma_i}(H_i)]\le [\hat \Gamma_i: \ker \phi]=|\op{im}(\phi)|\le |S(A_i)|\le |A_i|!$.

We now estimate $|A_i|$. Because $\pi_1(F_i)$ is normal in $\hat \Gamma_i$, $A_i$ is a subset of
all $\Bbb Z^k$-subgroups of $\pi_1(F_i)$. Let $\psi: \pi_1(F_i)\to \pi_1(F_i)/H_i$. Observe the following properties: (i) For any two subgroups of $H_1, H_2<\pi_1(F_i)$, $H_1=H_2$ if and only if $H_1\cap H_i=H_2\cap H_i$ and $\psi(H_1)=\psi(H_2)$.
(ii) $[H_i:H_1\cap H_i]\le [\pi_1(F_i):H_1]$. The two properties imply that $\pi_1(F_i)$ contains at most $c(a)$ many possible $\Bbb Z^k$-subgroups of index $a$, thus $|A_i|\le c(a)\le c(a(k))$, where $c(a)$ is a constant depending on $a$.

To see the claim, observe that the Gromov's short generators of $\hat \Gamma_i=\pi_1(\hat M_i)$ at $x_i$ are divided into two groups, one that their lengths are uniformly proportional to $\rho_i$ i.e., the short generators, $\alpha_i$ in $\pi_1(F_i)$, and $\gamma_i$ in $\hat \Gamma_i-\pi_1(F_i)$ such that
if $\hat \ell_i=\min\{|\gamma_i|,\,\, \gamma_i\in \hat \Gamma_i-\pi_1(F_i)\}$, then $\frac {\rho_i}{\hat \ell_i}\to 0$ as $i\to \infty$. Because $|\op{II}(N_i)|\le c(n)$, each $\gamma_i$-conjugation on $\alpha_i$ has a bounded distortion, the $\gamma_i$-conjugation preserves $\pi_1(F_i)$.
\qed\enddemo

\proclaim{Lemma 9} {\rm (Local product $T^k$-bundles)}  Let the assumptions be as in Lemma 8.

\noindent {\rm (9.1)} Let $\Bbb Z^k$ be as in {\rm (8.7)}. Then the centralizer of $\Bbb Z^k$,
$[\Gamma_i: C_{\Gamma_i}(\Bbb Z^k)]\le c(n)$.

\noindent {\rm (9.2)} For any $\delta>0$, there is $N(\delta)>0$, such that for $i\ge N(\delta)$, $x_i\in \hat M_i=
\tilde M_i/C_{\hat \Gamma_i}(\Bbb Z^k)$,
$$d_{\op{GH}}(B_{30}(x_i,(d\rho_i)^{-1}\hat M_i),B_{30}(T^k\times \Bbb R^{n-k},e\times 0))<\delta,$$
and $k$, $\op{diam}(T^k)\le 1$ and $\op{injrad}(T^k)\ge d^{-2}$ hold for all $x_i$.
\endproclaim

Note that (9.2) is required to construct a principal $T^k$-bundle on $\tilde M_i/C_{\Gamma_i}(\Bbb Z^k)$. Note that
in general a complete flat manifold $X$ with a soul flat $T^k$ may not isometric to the metric product,
$T^k\times \Bbb R^{n-k}$ (see Example 10).

\demo{Proof of Lemma 9}

(9.1) It suffices to prove (9.1) on $\hat M_i$; recall that by (8.7) we obtain that $N_{\hat \Gamma_i}(\Bbb Z^k)$
has a bounded index in $\hat \Gamma_i$. Fixing a canonical basis for $\Bbb Z^k$, $\gamma_{i,1},...,\gamma_{i,k}$ (consisting of closed geodesics in
$T^k$ at $x_i$), consider the holonomy representation induced by conjugation,
$\psi_i: N_{\hat \Gamma_i}(\Bbb Z^k)\to \op{Aut}(\Bbb Z^k)$.  We claim that the $\psi_i$-image is bounded (for all large $i$),
thus $\ker(\psi_i)=C_{\hat \Gamma_i}(\Bbb Z^k)$ has index in $\hat \Gamma_i$, $|\op{Im}(\psi_i)|$.

As a preparation, we first show that for any $\eta_i\in \hat \Gamma_i$ with $|\eta_i|\le 1$,
$$e^{-c(n)}\le \frac{d(\eta_i^{-1}\gamma_{i,j}\eta_i(\tilde x_i),\tilde x_i)}{d(\gamma_{i,j}(\tilde x_i),\tilde x_i)}\le e^{c(n)}.$$
We proceed by induction on $n$, starting with the trivial case $n=2$. For $n>2$,  we have
$N_i\to \hat M_i@>f_i>>\ell_iT^m$.  Applying induction on $N_i$, we assume that for $\alpha_i\in \Lambda_i$ with
$|\alpha_i|\le 1$, the above distortion estimate holds. Let $\alpha_i\beta_i$ denote a horizontal lifting of element from
$\pi_1(T^m)$ with $|\alpha_i\beta_i|\le 1$. Then
$$\split e^{-c(n)}\le&\frac{d((\alpha_i\beta_i)^{-1}\gamma_{i,j}(\alpha_i\beta_i)(\tilde x_i),\tilde x_i)}{d(\gamma_{i,j}(\tilde x_i),\tilde x_i)}\\
&=\frac{d(\alpha_i^{-1}\gamma_{i,j}\alpha_i(\beta_i(\tilde x_i)),\beta_i(\tilde x_i))}{d(\gamma_{i,j}(\beta_i(\tilde x_i)),\beta_i(\tilde x_i))}\cdot  \frac{d(\gamma_{i,j}(\beta_i(\tilde x_i)),\beta_i(\tilde x_i))}{d(\gamma_{i,j}(\tilde x_i),\tilde x_i)}\le e^{c(n)},\endsplit$$
because the distortion of $\alpha_i^{-1}\gamma_{i,j}\alpha_i$ (resp. $\gamma_{i,j}$) applies to every point in $\tilde N_i$ (see (8.5)), which contains $\tilde x_i$ and $\beta_i(\tilde x_i)$.

Let $\{\alpha_{i,j}\}_{j=1}^{a(n)}$ denote a Gromov's short generators for $N_{\hat \Gamma_i}(\Bbb Z^k)$ at $\tilde x_i$. Then $S_i=\{\psi_i(\alpha^\pm_{i,j})\}_{j=1}^{a(n)}$ is a symmetric set of generators for $\psi_i(N_{\hat \Gamma_i}(\Bbb Z^k))$, equipped with the word length in terms of $S_i$. Because $|\alpha_{i,j}|\to 0$ as $i\to \infty$, words with length
$\le 1$ looks like more and more, while the above estimate shows that the total number is  bounded above by $L\approx \frac{\op{vol}(\b B^k_{d^2c(n)+1}(0))}{\op{vol}(B_1^n(0))}$ (rescaling $|\gamma_{i,1}|=1$, $|\gamma_{i,k}|\le d^2$). In particular, $|(S_i)^L|<L$, where $(S_i)^L$ denotes the set of elements
of word length $\le L$. Consequently, $|\op{Im}(\psi_i)|<L$ for all large $i$ (\cite{KPT}).
Then
$$[\Gamma_i:C_{\Gamma_i}(\Bbb Z^k)]\le [\Gamma_i,C_{\hat \Gamma_i}(\Bbb Z^k)]\le [\Gamma_i,\hat \Gamma_i]
\cdot [\hat \Gamma_i: C_{\hat \Gamma_i}(\Bbb Z^k)]\le c(n).$$

(9.2) Arguing by contradiction, assuming $\epsilon>0$ and a sequence, $\hat M_i=\tilde M_i/C_{\hat \Gamma_i}(\Bbb Z^k)$, $x_i\in \hat M_i$, such that
$B_{30}(x_i(d\rho_i)^{-1}\hat M_i)$ is at least $\epsilon$-away from any $B_{30}(T^k\times \Bbb R^{n-k})$, with $\op{diam}(T^k)\le 1$ and
$\op{injrad}(T^k)\ge d^{-2}$. Passing to a subsequence, we may assume that $$(\rho_i^{-1}\hat M_i,x_i)@>\op{GH}>>
(X,x).$$

We shall show that $X=T^k\times \Bbb R^{n-k}$, a contradiction (see (8.6)). We proceed by induction on $n$, starting with
$n=2$. In this case, $\hat M_i=T^2$ (see Example 10). For $n>2$, as seen in the proof of (8.6), bases on
$N_i\to \hat M_i\to \ell_iT^m$ we have that $(\rho_i^{-1}\hat M_i, x_i)@>\op{GH}>>(Y\times \Bbb R^m,x)$ and
$(\rho_i^{-1}N_i,x_i)@>\op{GH}>>(Y,x)$.
Applying induction on $N_i$ ($\Lambda_i=C_{\Lambda_i}(\Bbb Z^k)$) we conclude that $Y=T^k\times \Bbb R^{n-m}$, thus
$X=T^k\times \Bbb R^{n-k}$.
\qed\enddemo

With the above preparation, we are ready for

\demo{Proof of Theorem 3'}

Without loss of generality, we may assume that all $M_i$ satisfies the additional regularities in Theorem 2, say $\ell=2$.  Fixing $i$ large, let
$\hat M:=\hat M_i=\tilde M_i/C_{\Gamma_i}(\Bbb Z^k)$.

(3.1)'  We divide the construction of a principal $T^k$-bundle in two steps: Step 1. Using (9.2), we construct a locally finite open
cover for $\hat M_i$, $\{U_\alpha\}$, $U_\alpha$ admits
an almost isometric $T^k$-action that is determined by (collapsed) metric, such that on $U_\alpha\cap U_\beta\ne\emptyset$,
the two $T^k$-actions are at least $C^1$-close.  Step 2. Glue $\{(U_\alpha,T^k)\}$ to a $T^k$-bundle with an affine structural
group. Because $\pi_1(\hat M_i)$ is the centralizer of $\pi_1(T^k)$, the structural group is trivial i.e., the $T^k$-bundle is a principal.

Step 1. Let $\{p_\alpha\}$ denote a $10$-net on $(d\rho_i)^{-1}\hat M$. Fixing $\epsilon>0$, by Lemma 9
($i$ large) we may assume a diffeomorphism, $f_\alpha: B_{30}(p_\alpha,(d\rho)^{-1}\hat M)\to B_{30}(0_\alpha,T^k\times \Bbb R^{n-k})$,
such  that the pullback flat metric is $C^3$ $\epsilon$-close. Let $U_\alpha:=f^{-1}_\alpha(T^k\times \b B_{10}^{n-k}(0_\alpha))$,
$V_\alpha:=f_\alpha^{-1}(T^k\times \b B_{12}^{n-k}(0_\alpha))$, where $\op{diam}(T^k)\le 1$ and $\op{injrad}(T^k)\ge d^{-2}$.
As discussed in the outline of the proof of Theorem 3', a choice of a base point $e\in T^k$ and a canonical
basis for $\pi_1(T^k,e)$ uniquely determine a group structure on $T^k$, hence $U_\alpha$ (resp. $V_\alpha$)
has a trivial principal $T^k$-bundle structure, $P_\alpha=\op{proj}_2\circ f_\alpha: U_\alpha\to \b B_{10}^{n-k}(0_\alpha)$.

For $U_\alpha\cap U_\beta \ne\emptyset$, let $W_\alpha=P_\alpha(V_\alpha\cap (U_\beta-U_\alpha))$, $W_\beta=B_{10}^{n-k}(0)-P_\beta(U_\alpha)$, and define a map,
$\phi_{\alpha\beta}: W_\alpha\to \b B_{10}^{n-k}(0_\beta)$,
$\phi_{\alpha\beta}(z)=\op{cm}(z)$, the center of mass of $P_\beta(P_\alpha^{-1}(z))$, $z\in W_\alpha$. It is clear that $\phi_{\alpha\beta}$
is an embedding, by which we glue together $\b B^{n-k}_{10}(0_\alpha)$ and $\b B_{10}^{n-k}(0_\beta)$: $\b B_{10}^{n-k}(0_\alpha)\cup_{W_\alpha\sim \phi_{\alpha\beta}(W_\alpha)} W_\beta$. Let $Z=P_\alpha^{-1}(W_\alpha)\cup P_\beta^{-1}(\phi_{\alpha\beta}(W_\alpha))$. Fixing a common base point $x_{\alpha\beta}\in U_\alpha\cap U_\beta$ and a canonical basis, up to an automorphism of $T^k$ we may view the two $T^k$-actions are from (a single) $T^k$, and
the two $T^k$-actions are at least $C^0$-close.

Using the fact that the pullback metrics on $U_\alpha$ and $U_\beta$ are $C^3$-close by which the two isometric $T^k$-actions
are isometric, by \cite{GK} the two $T^k$-actions are at least $C^1$-close.

\vskip2mm

Step 2. We will patch $\{(U_\alpha,T^k)\}$ together to obtain a $T^k$-bundle with affine structural group.  Observe that the desired result
follows from the gluing construction (Lemmas 1.4 and 1.5 of \cite{CG2}), based on the stability of compact group actions (\cite{GK}).
For the sake of a self-contained proof, here we present an elementary gluing construction based on the center
of mass technique in \cite{GK}.

For  each $z\in \phi_{\alpha\beta}(W_\alpha)$, we define a diffeomorphism, $\psi_z: P_\alpha^{-1}(z)\to P_\beta^{-1}(z)$, the
nearest point projection, which is smooth in $z$. We then replace $\psi_z(x)$ by the center of
mass of $\{t^{-1}\phi_z(tx),\, t\in T^k\}$, still denoted by $\psi_z(x)$, thus $\psi_z: P_\alpha^{-1}(z)\to P_\beta^{-1}(z)$ is also
$T^k$-conjugate. Let $\gamma_{x,\psi_z(x)}$ denote the unique normal minimal geodesic from $x$ to $\psi_z(x)$ (with respect to
a pullback flat metric on $B_{30}(x_{i,\alpha})$).

To patch $U_\alpha$ and $U_\beta$ together, we fix a transition function (which may be thought as the normalized distance function on $A_{11}^{n-k}(0_\alpha)=\b B_{11}^{n-1}(0_\alpha)-\b B_{10}^{n-k}(0_\alpha)$), on $\b B_{10}^{n-k}(0_\alpha)\cup_{W_\alpha\sim \phi_{\alpha\beta}(W_\alpha)} W_\beta$, as follows:
$$\chi(v)=\cases 0 & v\in \b B^{n-k}_{10}(0_\alpha) \\ 1 & v\in W_\beta-\phi_{\alpha\beta}(W_\alpha)\endcases,$$
We define $\Phi: [0,1]\times (U_\alpha\cup U_\beta)\to U_\alpha\cup U_\beta$:
$$\Phi(t,x)=\cases x & x\in U_\alpha\cup P_\beta^{-1}(W_\beta-\phi_{\alpha\beta}(W_\alpha))\\
\gamma_{x,\psi_{P_\beta(x)}(x)}
(\chi(P_\beta(x))t) & \text{otherwise} \endcases$$
Clearly, $\Psi(t,x)$ is an isotopy between $\op{id}_{U_\alpha\cup U_\beta}=\Phi(0,x)$ and $\Phi(1,x)$, and the new $T^k$-action on $U_\alpha\cup V_\alpha\cap U_\beta$ by
$\Phi(1,T^k)$, will agree with the $T^k_\beta$-action on the boundary of
$U_\alpha\cup (V_\alpha-U_\alpha)\cap U_\beta$, thus the two principal $T^k$-bundles
are glued into a $T^k$-principal $T^k$-bundle (a tool for a gluing in a general
situation is the isotopy extension theorem).

\vskip1mm

Because $\{U_\alpha\}$ is a finite set, we apply induction on the number of the gluing operations. By the above, we glue $U_1$ and
$U_2$ to form a $T^k$-bundle (with an affine structural group). Assume that we obtain a $T^k$-bundle on $U=U_1\cup\cdots \cup U_m$
(with an affine structural group). To glue $U_{m+1}$ to $U$,
it suffices to check that the gluing region in $W\cap U_{s+1}$ satisfies the $C^3$ $\epsilon$-closeness to
$B_{20}(T^k\times \Bbb R^{n-k})$. Because each gluing may produce an error on the $C^3$ $\epsilon$-closeness, say
$C^3$ $c(n)\epsilon$-close, a cumulation of errors could add up to excess $C\epsilon$ for any constant $C$.
This bad possibility is ruled out by the fact that in a region where multiple gluing operations occur, the number of gluing
operation is bounded above by $\le a(n)$.  Hence, taking $\epsilon$ suitably small we can apply the above model situation and
glue $W$ and $U_{s+1}$ to obtain a $T^k$-bundle with an affine structural group on $W\cup U_{m+1}$.

Finally, we explain that a desired $C^1$-close $T^k$-invariant metric is obtained by averaging the metric under the almost isometric
$T^k$-action: based on the additional regularities of the metric on $M$, the second fundamental form
and the $A$-tensor the principal $T^k$-bundle has bounded covariant derivative, thus the $T^k$-invariant metric satisfies the same regularity at least one order less than the original metric. By the O'Neil's
Riemannian submersion formula, it is clear that the quotient metric on $\hat M_i/T^k$ has bounded sectional curvature.

\vskip2mm

(3.2)'  The first part has been proved in (8.4). Consider $\ell_i^{-1}M_i@>\op{GH}>> \Bbb R^n/G$ in the proof of (8.1). Because $\Gamma_i$ is nilpotent, $T^\ell=\Bbb R^n/G$ (e.g., if the isometric $G$-action is not free, then $\Bbb R^n/G$ is not compact). Applying Theorem 7,
there is a fiber bundle map, $f_i: M_i\to T^{\ell}$, such that $|\op{II}(f_i)|\to 0$ with respect to $\ell_i^{-1}M_i$ (comparing with the proof of (8.2)). Because
for all $\tilde x_i$ in the fundamental domain at $\tilde p_i$, $\frac{d(\gamma_i(\tilde p_i),\tilde p_i)}{d(\gamma_i(\tilde x_i),\tilde x)}$ is scaling invariant,
this ratio is almost one is clear.
\qed\enddemo

We conclude the paper by supplying the following example mentioned several times through the
above proofs.

\example{Example 10} We construct one-parameter family of flat metrics, $g_\epsilon$, on a Klein bottle, such that
$(K^2,g_\epsilon)@>\op{GH}>>\op{pt}$, and the blow-up sequence, $(K^2,\ell_{\epsilon}^{-1}g_\epsilon)@>\op{GH}>>[0,1]$, not a manifold. Moreover, a pointed GH-limit, $(K^2,\rho_\epsilon^{-1}g_\epsilon,x)@>\op{GH}>>(X,x)$ is a complete flat surface which is either
isomorphic to $S^1\times \Bbb R^1$ or a twisted $\Bbb R^1$-bundle over $S^1$ (infinite M\"obius band), depending on $x$
(see (8.8) and Lemma 9).

First, $(K^2,g_0)$ admits an isometric $T^1$-action with two exceptional $T^1$-orbits, where $g_0$ is a flat metric.
Hence, $g_0$ at each point can be written $g_0=ds^2+(ds^2)^\perp$. For $\epsilon\in (0,1]$, define $g_\epsilon'=
\epsilon^2ds^2+(ds^2)^\perp$ and $g_\epsilon=\sqrt \epsilon g_\epsilon'$.  Then $(K^2,g_\epsilon)@>\op{GH}>>\op{pt}$.
It is clear that $\ell_\epsilon\approx \op{diam}(K^2,g_\epsilon)$, and
$$(K_2,\ell_\epsilon^{-1}g_\epsilon)@>\op{GH}>>[0,1],\quad (T^2,\ell_\epsilon^{-1}g_\epsilon)@>\op{GH}>>T^1.$$
Moreover, for any $x\in K^2$, $\op{injrad}(x)\approx \epsilon^{\frac 32}$, and
$$(K^2,\rho_\epsilon^{-1}g_\epsilon),x)@>\op{GH}>>\cases (S^1\times \Bbb R^1,x_\infty) & T^1(x) \text{ is regular} \\ (M_\infty,x_\infty) & T^1(x) \text{ is exceptional} \endcases,$$
where $M_\infty$ is a twisted $\Bbb R^1$-bundle over $S^1$ (an infinite M\"obius band).

\endexample

\noindent {\bf Acknowledgment}. The present paper is a reversion of ArXiv:1906.03377v1, in order to make the new proof self-contained and transparent to graduate students in Metric Riemannian geometry. The author would like to thank Jeff Cheeger for his suggestion to revise the early version for the above goal, and for his encouragement and support through many years.

\vskip4mm

\head Appendix.  Proof of Theorem 7
\endhead

\vskip4mm

A fiber bundle structure as in Theorem 7 has been generalized under various weak curvature conditions; however one may not have that the induced metric on a fiber is almost flat. Let's first briefly describe methods of constructions.

\vskip1mm

\noindent (A1) Smoothing method with $|\op{sec}_M|\le 1$ (\cite{CFG}, \cite{Fu1}). Given $M$ a metric $g_\epsilon$ with higher regularity that
is $C^1$ $\epsilon$-close to the original metric (Theorem 2), and a continuous $\epsilon$-GHA $h: M\to N$, using the method of center of mass
one smooths $h$ with respect to $g_\epsilon$ to be a desired bundle map, $h_\epsilon$ (called $\epsilon$-regular), satisfying (7.1)-(7.3).
Hence, $h_\epsilon$ satisfies (7.1)-(7.3) with respect to the original metric as well.

\vskip1mm

\noindent (A2) Imbedding method with $\op{sec}_M\ge -1$ (\cite{Fu4}, \cite{Ya}, \cite{Ro1,2}). Given a $\rho$-net, $\{\bar x_\alpha\}_{\alpha=1}^s\subset N$ and a cut-off function $h$, $\op{supp} (h)\subset [0,\rho)$, define a smooth map, $\Phi: N\to \Bbb R^s$, $\Phi(x)=(h\circ d_{\bar x_\alpha}(x))$, and a $C^1$-map, $\Psi: M\to \Bbb R^s$, $\Psi(x)=(h\circ \hat d_{x_\alpha}(x))$, where $\hat d$ is the average of $d_{x_\alpha}$ over  $B_{\frac \rho{10}}(x_\alpha)$, and $x_\alpha$ is a lift of $\bar x_\alpha$.  Then $\Phi: N\to \Bbb R^s$ is an embedding, and $f=\Phi^{-1}\circ \op{Proj}\circ \Psi: M\to N$ is a $C^1$-fiber bundle map satisfying (7.1) and a weak version of (7.2), where $\op{Proj}$ denotes the projection to a nearest point in $\Phi(N)$ (a fiber may not be an infra-nilmanifold; e.g., a sphere).

\vskip1mm

\noindent (A3) Gluing $(m,\delta)$-splitting maps with $\op{Ric}_M\ge -(n-1)$ and the universal cover of any $\rho$-ball on $M$
is non-collapsed (\cite{Hu}). Given an $\rho$-net, $\{\bar x_\alpha\}\subset N$, and let $x_\alpha$
be a lift of $\bar x_\alpha$. By \cite{CC}, one gets that a smooth $(m,\delta)$-splitting map, $H_\alpha: B_\rho(x_\alpha)\to B_\rho(0)\subset T_{\bar x}N$.
Let $\{f_\alpha\}$ be a partition of unity associate to the locally finite open cover, $\{B_\rho(x_\alpha)\}$, and the center of mass
method, one can glue $\{\exp\circ H_\alpha\}$ together to a smooth $H: M\to N$. Using the non-degeneracy of an $(n,\delta)$-splitting map (\cite{CJN}),
Huang showed that each $(m,\delta)$-splitting map $H_\alpha$ is non-degenerate and satisfies (7.1), and $H$ is also non-degenerated (the weak regularity fails to identify a fiber as an infra-nilmanifold).

\vskip1mm

\noindent (A4) Successively blowing up and gluing method (\cite{Ro3}) under the same conditions of (A3); a baby version is similar to the construction of the principal $T^k$-bundle in the proof of Theorem 3'.
The result is a fiber bundle map, $f: M\to N$, satisfying (7.1), such that a $f$-fiber diffeomorphic to an infra-nilmanifold, and the structural group.

\vskip2mm

In view of the above, our proof below does not rely on a smoothing method, nor involving in
an tedious calculus estimate. Using that $|\op{sec}_M|\le 1$, we will get the desired (local) regularities (for `free') based on the following lemma (see Step 2), via directly
verification of boundary conditions.

\proclaim{Lemma A5} Let $\Omega\subset \Bbb R^n$ be a bounded domain with smooth boundary,  let $g_i$ be two
Riemannian metrics on $\Omega$, and let $u_i$ be solution of the Dirichlet problems:
$$\cases \Delta_i u_i=0 & \Omega\\
u_i=f_i & \partial \bar \Omega \endcases,\qquad i=1, 2$$
Assume that for $k\ge 0$, $|g_1-g_2|_{C^{k,\alpha}(\Omega)}<\epsilon$ and $|f_1-f_2|_{C^{k+1,\alpha}(\partial \Omega)}<\epsilon$.
Then $|u_1-u_2|_{C^{k+1,\alpha}(\Omega)}\le \Psi(\epsilon)$.
\endproclaim

Observe that $v=u_1-u_2$ satisfies $\Delta _1v=(\Delta_2-\Delta_1)u_2$ on $\Omega$ and $v|_{\partial \Omega}=f_1-f_2$.
Hence Lemma A5 follows from a regularity estimate of an elliptic equation \footnote{this observation was pointed out to the
author by Zhengchao Han. \hfill{$\,$}}.


\demo{Proof of Theorem 7}

By a standard compactness argument, it suffices to consider a sequence of compact $n$-manifolds, $M_i@>\op{GH}>>N$, such that
$|\op{sec}_{M_i}|\le 1$, $|\op{sec}_N|\le 1$ and $\op{injrad}(N)\ge 2\rho$, and prove that for $i$ large there is a smooth map, $f_i: M_i\to N$, satisfying
(7.1)-(7.3) (note that (7.2) and (7.3) are local properties, and (7.2) implies that $f_i$ is a submersion, hence $f_i$ is a bundle map).


\vskip1mm

Step 1. Construct a smooth map, $f_i: M_i\to N$ (\cite{Hu}). Fixing $0<\delta\le 1$,
let $\{x_\beta\}$ be an $\delta\rho$-net on $N$. 
For each $x_\alpha$, let $e_1,...,e_m$ be a standard basis for $T_{x_\beta}N$, let $q_{\beta,j}=\exp_{x_\beta}\rho(e_j)$, $1\le j\le m$. For each $(i,\beta)$, let $(x_{i,\beta},\{q_{i,\beta,j}\}_{j=1}^m)\in M_i$ be a lift of $(x_\beta,\{q_{\beta,j}\}_{j=1}^m)$,  i.e. $x_{i,\beta}\to x_\beta$, $q_{i,\beta,j}\to q_{\beta,j}$, let $b_{i,\beta,j}=d(q_{i,\beta,j},\cdot)-d(x_{i,\beta},q_{i,\beta,j}): M_i\to\Bbb R$, and let $u_{i,\alpha,j}$ be the solution of the following Dirichlet problem (without loss of
the generality we may assume that $\partial \bar B_{\delta\rho}(x_{i,\beta})$ is smooth):
$$\cases \Delta u_{i,\beta,j}=0 &  B_{\delta\rho}(x_{i,\beta})\\ u_{i,\beta,j}=b_{i,\alpha,j} & \partial \bar B_{\delta\rho}(x_{i,\beta})\endcases
$$
We define a smooth map, $H_{i,\beta}: B_{\delta\rho}(x_{i,\beta})\to N$, $H_{i,\beta}(x_i)=\exp_{x_\beta}(\sum_{j=1}^m$
$u_{i,\beta,j}(x_i)e_j)$ (note that $H_{i,\beta}$ is a $(m,\delta)$-splitting map studied in \cite{CC}).
For $i$ large, $\{U_{i,\beta}=H_{i,\beta}^{-1}(B_{\frac{4\delta\rho}5}(x_\beta))\}$ is a locally finite open cover for $M_i$.

We will glue the family of local (fiber bundle) maps, $\{(U_{i,\beta},H_{i,\beta})\}$, by the method of center of mass: take a partition of unity, $\{f_\alpha\}$, associate to the open cover $\{B_{\frac{4\delta\rho}5}(x_\beta)\}$ (see (2.5). For any $x_i\in M_i$ close to $x_\beta\in N$, the energy functional,
$$E(x_i,y)=\frac 12\sum_\beta f_\beta(H_{i,\beta}(x_i))d_N^2(H_{i,\beta}(x_i),y): \, B_\rho(x_\beta) \to \Bbb R.$$
is a finite sum of strictly convex functions in $B_\rho(x_\beta)$, $E(x_i,y)$ is strictly convex, thus $E(x_i,y)$
achieves a unique minimum in $\bar B_\rho(x_\beta)$, denote by $\op{cm}(x_i)$ (called the center of mass).
We define $f_i: M_i\to N$, $f_i(x_i)=\op{cm}(x_i)$.

\vskip1mm

Step 2. We will verify that $H_{i,\beta}$ satisfies (7.1) and (7.3) with $\delta=\frac 14$, and
$\delta=\delta(\epsilon)$ suitably small and fixed, $\delta^{-1}M_i\to \delta^{-1}N$, $H_{i,\beta}$ satisfies (7.2) which is scaling invariant. Note that (7.2) and (7.3) are local properties, so our
strategy is to verify (7.2) and (7.3) to the lifting of $H_{i,\beta}$ on a
`covering space'.

Let $B_{\pi}(0_{i,\beta})\subset T_{x_{i,\beta}}M_i$ equipped with the pullback metric by
$\exp_{x_{i,\beta}}$, $g_{i,\beta}^*$, let $\hat u_{i,\beta,j}(\tilde x_i)=u_{i,\beta,j}(\exp_{x_{i,\beta}}(\tilde x_i))$, $\hat b_{i,\beta,j}(\tilde x_i)=b_{i,\beta,j}(\exp_{x_{i,\beta}}(\tilde x_i))$, defined on $V_{i,\beta}=\bigcup_{\gamma_i\in \Lambda_{i,\beta}}\gamma_i(D_{i,\beta}\cap B_{\delta\rho}(0_{i,\beta}))$, $B_{\delta\rho}(0_{i,\beta})\subset V_{i,\beta}\subset B_{2\delta\rho}(0_{i,\beta})$, where $D_{i,\beta}\subset B_{\frac \pi2}(0_{i,\beta})$ is a `fundamental domain' of $\exp_{x_{i,\beta}}$ at $0_{i,\beta}$ (i.e., $\exp_{x_{i,\beta}}: D_{i,\beta}\to B_{\frac \pi2}(x_{i,\beta})$
is bijection), $\Lambda_i\subset \pi_1(B_\epsilon(x_{i,\beta}))$ ($\epsilon<<\delta\rho$)
is a finite subset (roughly, one may imagine $V_{i,\beta}$ as a finite solid cylinder
with `boundary', $\partial V_{i,\beta}=V_{i,\beta}\cap \exp^{-1}_{x_{i,\beta}}(\partial B_{\delta\rho}(x_{i,\beta}))$, i.e., without the top/bottom caps). Then
$\hat u_{i,\beta,j}$ satisfies the following:
$$\cases \Delta \hat u_{i,\beta,j}=0 & V_{i,\beta} \\ \hat u_{i,\beta,j}=\hat b_{i,\beta,j} & 
\partial V_{i,\beta},\endcases$$
and $H_{i,\beta}$ satisfies (7.2) if and only if $\hat H_{i,\beta}(x)=(\hat u_{i,\beta,1}(x),...,\hat u_{i,\beta,m}(x))$ does, and
$\hat H_{i,\beta}$ is a `periodic' function with a `period',
$\hat H_{i,\beta} |_{D_{i,\beta}\cap B_{\delta\rho}(0_{i,\beta})}$. 

Restricting to $D_{i,\beta}$ (with $\delta=\frac 14$), we may assume that $\hat b_{i,\beta,j}=d(\rho e_j,\cdot)-\rho$. By standard Schauder interior estimate, (7.3) follows from that $|\hat u_{i,\beta_j}|_{C^\alpha(V_{i,\beta})}\le c(n)$. Note that $\hat u_{i,\beta,j}+\rho$ is positive harmonic on $V_{i,\beta}$, by the Cheng-Yau gradient estimate
$|\nabla \hat u_{i,\beta,j}|_{C^1(V_{i,\beta}\cap B_{\frac 12\delta\rho}(x_{i,\beta}))}\le c(n)$, if $|\hat u_{i,\beta,j}+\rho|_{C^0(D_{i,\beta}\cap B_{\frac 12\delta\rho}(0_{i,\beta}))}\le c(n)$. By the maximal principle $|\hat u_{i,\beta,j}+\rho|$ achieves the maximum on $\partial V_{i,\beta}$; thus the desired
bound is obvious, $|\hat b_{i,\beta,j}+\rho|_{C^0(\partial V_{i,\beta})}\le c(n,\rho)$.


In proving that $H_{i,\beta}$ satisfies (7.2), our approach is to construct a model
solution with respect to the Euclidean metric on $\Bbb R^n=T_{x_{i,\beta}}M_i$, and applying
Lemma A5.

Let $\b b_j(x)=d(\Bbb R(\rho e_j)^\perp,x)-\rho$ the distance function
to the hyperplane at $\rho e_j$, where $\Bbb R^\perp(\rho e_j)$ is the hyperplane
at $\rho e_j$ with the normal vector $e_j$. Then $\b b_j$ satisfies
$$\cases \Delta \b b_j=0 & V_{i,\beta} \\ \b b_j=\b b_j & \partial V_{i,\beta},\endcases$$
and $\b H_\beta(x)=(\b b_1(x),...,\b b_m(x)): V_{i,\beta} \,(\subset \Bbb R^n)\to \Bbb R^m$, is the projection; so (7.2) holds for any $\epsilon>0$. We claim that $\delta=\delta(\epsilon)$ is chosen small, $\b H_\beta$ is $C^{1,\alpha}$-close to $\hat H_{i,\beta}$ with respect to $\delta^{-1}N$ and $\delta^{-1}M_i$, thus $\hat H_{i,\beta}$ satisfies (7.2).
By Lemma A5, the claim follows from that restricting to $\partial V_{i,\beta}$,
$\hat b_{i,\beta,j}$ is $C^{1,\alpha}$-close to $\b b_j$.

For $\delta$ suitably small, the pullback metric $g_i^*$ on $B_{\delta\rho}(0_{i,\beta})$ is $C^{1,\alpha}$-close to the Euclidean metric (see the proof of (8.7)), thus it suffices to
check that $\hat b_{i,\beta,j}(x)=d_{A_{i,\beta}}(x)-\rho$ is $C^0$ $\epsilon$-close to
$\b b_j$, which equals to, restricting to $V_{i,\beta}$, the Busemann function along the $x_j$-axis; clearly, $\hat b_{i,\beta,j}$ approximates the Busemann function along the $x_j$-axis, for $\delta(\epsilon)<<1$ and $i$ large.

\vskip1mm

Step 3. By the construction of $H_i: M_i\to N$ and Step 2, it is clear that $H_i$ satisfies (7.1) and (7.3). It remains to check that $H_i$ satisfies (7.2).

Because (7.2) is local property, and because $H_{i,\beta}$ and $\b b_j$ are at least $C^{1,\alpha}$-close, it reduces to check (7.2) for the gluing of $\{\b H_\beta\}$ via $\{f_\beta\}$,
which is easily seen via a differentiation of implicit functions.
\qed\enddemo

\remark{Remark \rm A6} Inspecting the above proof, it is clear that if one assumes that
$g_i$ satisfies a little higher regularity, $|\nabla \op{Rm}(g_i)|\le C$, then (for $\delta$ suitably small; e.g. replacing $M_i$ with $\sqrt{\ell_i}M_i@>\op{GH}>> \op{pt}$) $\hat u_{i,\beta,j}$ is $C^{2,\alpha}$ close to $\b b_j$.
\endremark

\enddocument